\documentclass[12pt,a4paper]{article}

\usepackage{amssymb,amscd,amsmath,amstext,amsfonts}
\usepackage{a4wide}
\usepackage{theorem}
\usepackage{enumerate}

\newtheorem{thm}{Theorem}[section]
\newtheorem{defn}[thm]{Definition}

\newtheorem{prop}[thm]{Proposition}

\newtheorem{rem}[thm]{Remark}

\numberwithin{equation}{section}

\def\beginpf{\noindent {\bf Proof :} \quad}
\def\endpf{\rightline{$\square$}}

\def\N{{\mathbb N}}
\def\R{{\mathbb R}}
\def\T{{\mathbb T}}
\def\Z{{\mathbb Z}}

\newcommand{\ds}{\displaystyle}
\newcommand{\E}{\mathrm{E}}

\begin{document}

\title{Hyperinvariant subspace for weighted composition operator on $L^p([0,1]^d)$}
\author{George Androulakis \and Antoine Flattot \thanks{Part of this work was completed 
while the authors were supported by the 2008 SUMIRFAS conference.}}

\maketitle

\begin{abstract}
The main result of this paper is the existence of a hyperinvariant subspace of weighted 
composition operator $Tf=vf\circ\tau$ on $L^p([0,1]^d)$, ($1 \leq p \leq \infty$) when 
the weight $v$ is in the class of ``generalized polynomials'' and the composition map 
is a bijective ergodic transform satisfying a given discrepancy. The work is based on 
the construction of a functional calculus initiated by Wermer and generalized by Davie. 
\end{abstract}

\noindent
{\em  AMS 2000 subject classification}:
Primary- 47A15 ; secondary- 47A10, 47A60.\\
{\em Address}: Department of Mathematics, University of South Carolina, Columbia, SC 29208\\
{\em E-mails}: giorgis@math.sc.edu, flattot@math.sc.edu

\section{Introduction}
We study the existence of invariant (and even hyperinvariant) subspace for weighted composition operator on $L^p([0,1]^d)$, $1\leq p\leq\infty$, that means operator of the form
\[ Tf(x)=v(x)f(\tau(x)) ,\]
with $v\in L^\infty([0,1]^d)$ and $\tau:[0,1]^d\rightarrow[0,1]^d$. The aim is to give conditions on the weight $v$ and the composition application $\tau$ to obtain the exitence of such subspaces.

What we can first say is that if $\tau$ is not an ergodic transformation then there exists an invariant subspace. Indeed, if $\tau$ is not ergodic,  then there exists a Borel set, say $\Omega$, such that $| \Omega |>0$, $|[0,1]^d \setminus \Omega | >0$ and such that $\tau(\Omega)\subset \Omega$. Hence $T{\cal M}\subset {\cal M}$ where ${\cal M} = \chi_{\Omega}\,L^p([0,1]^d)$. Recall that given a closed subset $E$ of $[0,1]^d$, the subspace defined by $\{f\in L^p([0,1]^d):f_{|E}=0\}$ is called a spectral subspace. So, in other words, if $\tau$ is not ergodic, then $T$ has a spectral nontrivial closed invariant subspace.\\
Therefore, from now on, we shall consider ergodic map $\tau$ for the composition. In particular, $\tau$ is measure preserving. We will also assume that $\tau$ is a bijection.

Also, if $\inf_{x\in[0,1)}|v(x)|>0$ then the operator $T$ is invertible, and then the techniques of Wermer's theorem \cite{Wermer54} apply. Such work was done by Bletcher and Davie in \cite{BD90} when $\tau$ is an irrational rotation $\tau(x)=\{x+\alpha\}$ for certain irrational numbers $\alpha$ and where the weight $v$ does not vanish on $[0,1)$ and its moduli of continuity satisfies a given condition. MacDonald generalized this result in \cite{GMD96}, for a larger class of irrationals.

We will use also along this article the discrepancy associated to an ergodic map. Here we recall some standard terminology and facts, further details can be found in \cite{DT97,KN74}. Let consider, for $\omega=(x_n)$ a sequence of real numbers in $I=[0,1)$ and $E\subset I$, the set
\[ A(N,E,\omega)=\sharp \{ n\,;\, 1\leq n\leq N \,,\, x_n \in E \} .\]
Then the discrepancy associated to $\omega$ is defined as
\begin{equation} \label{DN}
D_N(\omega)=\underset{0\leq\alpha<\beta\leq1}{\sup}\left|\frac{A(N,[\alpha,\beta),\omega)}{N}-(\beta-\alpha)\right| .
\end{equation}
For a fixed $x\in [0,1)$ let $\omega(x)=(\tau^n(x))_{ n \geq 1}$, we denote by $D_N = \sup_{x \in [0,1)}D_N(\omega(x))$. Also, if $\tau$ is bijective, let $\omega'(x)= (\tau^{-n}(x))_{n \geq 1}$ and notice that $D_N(\omega'(x)) =D_N(\omega(\tau^{-N-1}(x))) \leq D_N$. The Birkhoff ergodic theorem states
that if $\tau$ is ergodic then for almost all $x$, $D_N(\omega(x)) \to 0$ as $N \to \infty$. Also it is known that for certain ergodic transformations
$\tau$, the function $x \mapsto D_N(\omega(x))$ is almost constant, i.e. $\limsup_{N \to \infty} \sup_x D_N(\omega(x))/ \inf_x D_N(\omega(x))$ $< \infty$. This is known to be the case for the irrational rotation. In this case, $D_N \to 0$ as $N \to \infty$.

The first part of the paper concerns our main theorem, that gives the existence of a hyperinvariant subspace for the weighted composition on $L^p([0,1])$, assuming a hypothesis on the discrepancy $D_N$ of the ergodic transform $\tau$. We then provide sufficient number theoretic conditions which guarantee that the discrepancy $D_N$ of an irrational rotation satisfies the assumptions of our main theorem. The last part deals with generalization of this work to $L^p([0,1]^d)$.

\section{Case one dimentional}

We introduce the class $\mathcal{P} \subset L^\infty([0,1])$ of functions on $[0,1]$ saying that $v\in \mathcal{P}$ if there exist some positive constants $s_0,\hdots,s_l$, a constant $C$ and map $\sigma_i$ such that
\begin{equation} \label{P} 
v(x)=C\prod_{i=0}^l(\sigma_i(x-x_i))^{s_i} ,
\end{equation}
where $x_i \in [0,1]$ and $\sigma_i(x)=x$ or $\sigma_i(x)=|x|$ for $i=0\hdots l$. Note that the absolute values may be needed for certain powers $s_i$ in order $v$ to be defined on $[0,1]$. We call this class the ``generalized polynomials'' (since we allow non-integer powers).

The main theorem of this paper is the following:
\begin{thm} \label{main}
Let $\tau$ be a bijective ergodic transformation of $[0,1]$ and $v\in \mathcal{P}$, we consider the operator on $L^p([0,1])$, $1\leq p\leq \infty$, defined by
\[ Tf(x)=v(x)f(\tau(x)).\]
Then, if the discrepancy $D_n$ of the sequence $(\tau^n)_{n\geq1}$ satisfies
\begin{equation}\label{cond_disc}
D_n =O\left(\frac{1}{\ln^{3+\epsilon} n}\right) 
\end{equation}
for some $\epsilon >0$, this operator has a hyperinvariant subspace.
\end{thm}

\beginpf 
We want to construct a functional calculus based on regular Beurling algebra $A_w$, associated to a subadditive weight $w$, as first appeared in Wermer \cite{Wermer54} and then generalized by Davie in \cite{Davie74} . This algebra is defined as follows:
\[ A_w =\left\{ \phi \in \mathcal{C}( \T ) \;;\; \sum_{n\in \Z}|\widehat{\phi}(n)|e^{w(|n|)} <\infty \right\} .\]
We refer to \cite[chap.5 \S 2]{CF68} for details about this Banach algebra. The important property for us is that if the weight satisfies $\sum_{n\geq1}w(n)/n^2<\infty$ then the algebra $A_w$ is regular.\\
For $\phi\in A_w$ we define an ``operator'' as follows:
\begin{equation} \label{series}
\phi(T)f=\sum_{n\in\Z} \hat{\phi}(n)T^nf .
\end{equation}
We need to give sense to this definition. The operator $T$ is not an invertible operator on $L^p([0,1])$, but since $\tau$ is a bijection and $v$ only vanishes in finitely many points, the inverse map $T^{-1}$ (as well as the other negative powers of $T$) of $T$ is well defined on the set of measurable functions. These are the maps that we will refer to when we write negative powers of $T$. Also since this series contains both positive and negative powers of $T$, if we want to obtain convergence of this series in some sense, then we need to ``normalize'' $T$ appropriately.  The correct ``normalization'' is to divide $T$ by its spectral radius. Of course, this does not effect the validity of Theorem~\ref{main}. For $v\in \mathcal{P}$, we can compute the spectral radius of the associated operator acting on $L^p([0,1])$ ($1 \leq p \leq \infty$), using \cite[Proposition~1.3]{GMD90} (this result was proved for $p=2$ but a glance at the formula of $T^n$ shows that the same result is true for $1 \leq p \leq \infty$):
\[ \begin{array}{rcl}
   r(T)&=& C \exp\left(\int_0^1 \ln\left(\prod_{i=0}^l\left|x-x_i\right|^{s_i}\right)dx\right) \\
       &=& C \exp\left( \sum_{i=0}^l -s_i(X_i+X'_i)\right) ,
   \end{array} \]
where, for $i=0\hdots l$
\begin{equation} \label{X}
X_i=x_i-x_i\ln x_i \quad \mbox{and} \quad X'_i=(1-x_i)-(1-x_i)\ln(1-x_i) ,
\end{equation}
with the convention that if $x_0=0$ or $1$ then $X_0+X'_0=1$. Note also that $1\leq X_i+X'_i \leq 1+\ln 2$ for all $i\in \{0,\hdots,l\}$, and that if $x_1=1-x_0$ then $X_0+X'_0=X_1+X'_1$.\\
Thus from now on, we assume that the weight $v$ in the formula of $T$ in Theorem~\ref{main} has the form $v(x)=1/C\exp\left( \sum_{i=0}^l s_i (X_i +X_i') \right) \prod_{i=0}^l ( \sigma_i (x-x_i))^{s_i}$ where $\sigma_i (x)=x$ or $\sigma_i (x) = | x |$.\\[0.2cm]
\indent In order to prove that the series (\ref{series}) converges in some sense, we first need to compute estimates on the powers $T^n$ for $n \in \Z$.\\
\noindent \textbf{Estimation of the bounds of $T^n$, $n\in\Z$:}\\
We prove the following theorem
\begin{thm} \label{bounds}
There exists a subadditive weight $w$ satisfying $\sum_{n\geq1}w(n)/n^2<\infty$ such that for $n\in\N$
\[ \|T^nf\|_{L^p([0,1])} \leq e^{w(n)}\|f\|_{L^p([0,1])} .\]
Moreover, for $t\in(0,1)$ there exists $n(t)\in \N$ and a set $E_t \subset [0,1]$ with $E_{t_1}\subset E_{t_2}$ for $t_1>t_2$ and $|E_t| \nearrow 1$ as $t\to 0$ such that 
\[ |T^{-n}f(x)| \leq \frac{1}{L_{n,x}}|f(\tau^{-n}x)| \quad \mbox{for all }x \in [0,1] \mbox{ and } n \in \N \]
where
\[ \sup_{x \in E_t} \frac{1}{L_{n,x}} < \infty \mbox{ for } n \in \N \mbox{ and } \sup_{x \in E_t} \frac{1}{L_{n,x}} \leq  e^{w(n)} \mbox{ for } \ n\geq n(t) .\]
\end{thm}

\beginpf
Let us first consider the particular case where the weight vanishes at one point, say $x_0\in (0,1)$ and the power
is equal to $1$, i.e.  we work with the weight
\[ v(x)=e^{X_0+X'_0}(x-x_0) .\]
We will be able to give an upper bound for $T^nf$ and $n\in\Z$.\\
An easy computation provide the following expressions of $T^n$:
\[ \begin{array}{l}
T^nf(x)=e^{n(X_0+X'_0)}\prod_{k=0}^{n-1}(\tau^k(x)-x_0) f(\tau^n(x)) \quad n\geq 1,\\
\ds{T^{-n}f(x)=\frac{1}{e^{n(X_0+X'_0)}\prod_{k=1}^{n}(\tau^{-k}(x)-x_0)} f(\tau^{-n}(x))} \quad n\geq 1,
\end{array} \]
and so we want to bound from above the quantity
\[ U_{n,x}=e^{n(X_0+X'_0)}\prod_{k=0}^{n-1}\left|\tau^k(x)-x_0\right| ,\]
and bound from bellow (with a lower bound not zero!) the following
\[ L_{n,x}=e^{n(X_0+X'_0)}\prod_{k=1}^{n}\left|\tau^{-k}(x)-x_0\right| .\]

The easiest case is the one of the positive powers of $T$ because we do not need to take care where $v$ vanishes.\\
\textit{The upper bound of the positive powers of $T$:}\\
We choose $\delta_{1,n}, \delta_{2,n} >0$ and $k_{1,n}=x_0/\delta_{1,n}, k_{2,n}=(1-x_0)/\delta_{2,n} \in\N$ to obtain a partition of the intervals $I_1=[0,x_0)$ and $I_2=[x_0,1)$. Let denote $\Delta_n=\max(\delta_{1,n},\delta_{2,n})$ and $\delta_n=\min(\delta_{1,n},\delta_{2,n})$ and assume that
\begin{equation}\label{upper_bound_Delta_n/delta_n}
\frac{\Delta_n}{\delta_n}\leq 1+\frac{1}{\ln^2 n} .
\end{equation}
Let, for $i=1\hdots k_{1,n}$ and $j=1\hdots k_{j,n}$,
\[ \begin{array}{rcl} 
  n_i^1(x)&:=& \sharp \left\{ 1\leq k\leq n \,;\, \tau^{k}(x) \in [(k_{1,n}-i)\delta_{1,n},(k_{1,n}-i+1)\delta_{1,n}) \right\} \\
          &=& \sharp \left\{ 1\leq k\leq n \,;\, x_0-\tau^{k}(x) \in [(i-1)\delta_{1,n},i\delta_{1,n}) \right\} ,\\
  \mbox{and} \quad n_j^2(x)&:=& \sharp \left\{ 1\leq k\leq n \,;\, \tau^{k}(x) \in [x_0+(j-1)\delta_{2,n},x_0+j\delta_{2,n}) \right\} .
  \end{array} \]
We have
\[ \begin{array}{rcl}
     U_{n,x} &=& e^{n(X_0+X'_0)} \prod_{k ; \tau^{k}(x) \in I_1}(x_0-\tau^{k}(x))\prod_{k ; \tau^{k}(x) \in I_2}(\tau^{k}(x)-x_0) \\ 
          &\leq& e^{n(X_0+X'_0)} \prod_{i=1}^{k_{1,n}}(i\delta_{1,n})^{n_i^1(x)} \prod_{j=1}^{k_{2,n}}(j\delta_{2,n})^{n_j^2(x)} \\
          &=& e^{n(X_0+X'_0)} \delta_{1,n}^{n_1^1(x)+\cdots+n_{k_{1,n}}^1(x)}\delta_{2,n}^{n_1^2(x)+\cdots+n_{k_{2,n}}^2(x)} \prod_{k=2}^{k_{1,n}}k^{n_k^1(x)}\prod_{k=2}^{k_{2,n}}k^{n_k^2(x)} \\
          &\leq& e^{n(X_0+X'_0)} \Delta_n^{n} \prod_{k=2}^{k_{1,n}}k^{n_k^1(x)}\prod_{k=2}^{k_{2,n}}k^{n_k^2(x)} .
    \end{array} \]
By definition of $D_n$ we get, for $i=1,\hdots,k_{1,n}$ and $j=1,\hdots,k_{2,n}$
\[ \frac{n_i^1(x)}{n}-\delta_{1,n} \leq D_n \quad \mbox{and} \quad \frac{n_j^2(x)}{n}-\delta_{2,n} \leq D_n ,\]
and $D_n \xrightarrow[n\to\infty]{}0$.\\
Thus,
\[ U_{n,x} \leq e^{n(X_0+X'_0)} \delta_n^{n} (k_{1,n} !)^{n(D_n+\delta_{1,n})} (k_{2,n} !)^{n(D_n+\delta_{2,n})}=U_n .\]
Let us look first at $(k_{1,n} !)^{n(D_n+\delta_{1,n})}$.\\
We use the asymptotical development of the gamma function:
\begin{equation}\label{dev_gamma}
\Gamma(x)=\left(\frac{x}{e}\right)^x \sqrt{2\pi x}\left(1+\frac{1}{12x}+\frac{1}{288x^2}-\frac{139}{51840 x^3}+o(1/x^3)\right) .
\end{equation}
In the particular case of the factorial we obtain:
\[ n! \leq \left(\frac{n}{e}\right)^n\sqrt{2\pi n}\left(1+\frac{1}{12n}+\frac{1}{288n^2}\right) \leq e\left(\frac{n}{e}\right)^n\sqrt{2\pi n} .\]
Hence, using this and the fact that $\delta_{1,n} \rightarrow 0$, we get for $n$ large enough:
\[ \begin{array}{rcl}
    (k_{1,n} !)^{n(D_n+\delta_{1,n})} &\leq& \exp\left(n(D_n+\delta_{1,n}) \ln\left(e\left(\frac{k_{1,n}}{e}\right)^{k_{1,n}}\sqrt{2\pi k_{1,n}}\right)\right) \\
      &=& \exp\left(n(D_n+\delta_{1,n}) \left[\frac{x_0}{\delta_{1,n}}\ln\left(\frac{x_0}{e\delta_{1,n}}\right) + 1/2\ln\left(\frac{2\pi x_0}{\delta_{1,n}}\right) + 1\right]\right) \\
      &=& \exp\left(n(D_n+\delta_{1,n}) \left[-\frac{X_0}{\delta_{1,n}} -\frac{x_0\ln\delta_{1,n}}{\delta_{1,n}}+1/2\ln(2\pi x_0)-1/2\ln\delta_{1,n}+1\right]\right) \\
      &=& \exp\left( -nX_0-nx_0\ln \delta_{1,n}-\frac{nD_n}{\delta_{1,n}}(1+x_0\ln\delta_{1,n})\right. \\
        && \hspace*{4cm}    \left. +n(D_n+\delta_{1,n}) \left[1/2\ln(2\pi x_0)-1/2\ln\delta_{1,n} +1\right]\right) \\
      &\leq& \exp\left( -nX_0-nx_0\ln \delta_{1,n}-\frac{nD_n}{\delta_{1,n}}(x_0\ln\delta_{1,n})+n(D_n+\delta_{1,n}) \left[-\ln\delta_{1,n}\right]\right) \\
      &=& \exp\left( -nX_0-nx_0\ln \delta_{1,n}+nD_n\left[-\frac{x_0\ln\delta_{1,n}}{\delta_{1,n}}-\ln\delta_{1,n}\right]-n\delta_{1,n}\ln\delta_{1,n} \right) \\
      &\leq& \exp\left( -nX_0-nx_0\ln \delta_{1,n}-2nD_n\frac{x_0\ln\delta_{1,n}}{\delta_{1,n}}-n\delta_{1,n}\ln\delta_{1,n} \right) \\
      &\leq& \exp\left( -nX_0-nx_0\ln \delta_n-2nD_n\frac{x_0\ln\delta_n}{\delta_n}-2n\delta_n\ln\delta_n \right).   
   \end{array}\]
Doing the same with $k_{2,n}$ we obtain:
\[ (k_{2,n} !)^{n(D_n+\delta_{2,n})} \leq \exp\left( -nX'_0-n(1-x_0)\ln \delta_n-2nD_n\frac{(1-x_0)\ln\delta_n}{\delta_n}-2n\delta_n\ln\delta_n \right). \]
Therefore it comes
\[ \begin{array}{rcl}
 U_n &\leq& \exp\left(n\ln\frac{\Delta_n}{\delta_n}-2nD_n\frac{\ln\delta_n}{\delta_n}-4n\delta_n\ln\delta_n\right)\\
     &\leq& \exp\left(\frac{n}{\ln^2 n}-2nD_n\frac{\ln\delta_n}{\delta_n}-4n\delta_n\ln\delta_n\right),
    \end{array} \]
since, by (\ref{upper_bound_Delta_n/delta_n}), $\ln(\Delta_n/\delta_n)\leq \ln(1+1/\ln^2n)\leq 1/\ln^2n$.\\
We choose $\delta_n$ such that $n\delta_n\ln\delta_n/n^2$ is summable, a good choice is $\delta_n=1/\ln^{1+\epsilon_1}n$. For this value of $\delta_n$ we set
\[ \widetilde{w_1}(n)=2nD_n\ln^{1+\epsilon_1}n \ln(\ln^{1+\epsilon_1}n)+2n\frac{\ln(\ln^{1+\epsilon_1}n)}{\ln^{1+\epsilon_1}n} .\]
Hence, if
\[ D_n =O\left(\frac{1}{\ln^{2+\epsilon}n}\right) ,\]
where $\epsilon$ is an arbitrary positive number, we get that $\widetilde{w_1}/n^2$ is summable (we can change the $\epsilon_1$ if necessary).\\
Fix $\epsilon>0$ and take $\epsilon_1<\epsilon$, for the choice of $D_n=1/\ln^{3+\epsilon}n$ (the choice of the power $3+\epsilon$ will be clearer in the light of the lower bound) we obtain:
\[ \widetilde{w_1}(n)=2n(1+\epsilon_1)\frac{1}{\ln^{2+\epsilon-\epsilon_1}n} \ln(\ln n)+2n\frac{\ln(\ln^{1+\epsilon_1}n)}{\ln^{1+\epsilon_1}n} .\]
Moreover, for any $\eta >0$, we have for $x$ large enough (depending on $\eta$) :
\[ \frac{\ln x}{x} \leq \frac{1}{x^{1-\eta}} .\]
So with a good choice of $\epsilon_1$ we can obtain:
\[ \widetilde{w_1}(n) \leq 4\frac{n}{\ln^{2+\epsilon/2}n}+ 2\frac{n}{\ln^{1+\epsilon/2}n} \leq 6\frac{n}{\ln^{1+\epsilon/2}n} .\]
Let 
\begin{equation} \label{wtilde}
\widetilde{w}(n)=15\frac{n}{\ln^{1+\epsilon/2}n} ,
\end{equation}
hence we get:
\[ U_n \leq \exp(\widetilde{w}(n)) .\]

Let us now see how to obtain a non zero lower bound.\\
\textit{The upper bound of the negative powers of $T$:}\\
We need to remove a little interval around $x_0$. For a fix $t>0$ and $n\in\N$ we work on the two subintervals $I_1=[0,x_0-t/n^3)$ and $I_2=[x_0+t/n^3,1)$.\\
We define the set
\begin{equation}\label{setEnxt}
E_{n,x_0,t}=\{ x\in[0,1] \,;\, |\tau^{-k}(x)-x_0| \geq t/n^3 \; \forall k=1,\hdots,n \} , 
\end{equation}
and $E_{x_0,t}=\bigcap_{n\geq 1} E_{n,x_0,t}$. Note that $(E_{x_0,t})_{t\geq0}$ is decreasing.\\
We have $E_{n,x_0,t}=\cap_{k=1}^n E_{n,x_0,t,k}$ with $E_{n,x_0,t,k}=\{ x\in[0,1] \,;\, |\tau^{-k}(x)-x_0| \geq t/n^3 \}$. Then, using the measure preserving of $\tau$, we get: 
\[ \begin{array}{rcl}
  | E_{n,x_0,t}|&=& 1-|E_{n,x_0,t}^c|\\
              &=& 1-|\cup_{k=1}^{n} {E_{n,x_0,t,k}^c}|\\
              &\geq& 1-t\sum_{k=1}^n 2/n^3 =1-2t/n^2
   \end{array} \]           
So $|E_{x_0,t}| \geq 1-\sum_{n\geq 1}2t/n^2=1-t\pi^2/3 \xrightarrow[t\to 0]{} 1$.\\
We work on $E_{x_0,t}$.\\
We define $\delta_{1,n}, \delta_{2,n} >0$ and $k_{1,n}, k_{2,n} \in\N$ such that
\[ k_{1,n}=\frac{x_0-t/n^3}{\delta_{1,n}} \quad \mbox{and} \quad k_{2,n}=\frac{1-(x_0+t/n^3)}{\delta_{2,n}} .\]
We partition the interval $I_1$ with the step size $\delta_{1,n}$, and $I_2$ with the step size $\delta_{2,n}$ . Once again we denote $\Delta_n=\max(\delta_{1,n},\delta_{2,n})$ and $\delta_n=\min(\delta_{1,n},\delta_{2,n})$ but we assume now that
\begin{equation}\label{upper_bound_Delta_n}
\Delta_n \leq \delta_n+\frac{1}{24\ln^3 n} .
\end{equation}
For $i=1,\hdots,k_{1,n}$ and $j=1,\hdots,k_{2,n}$ let
\[ \begin{array}{rcl}
    m_i^1(x)&:=& \sharp \left\{ 1\leq k\leq n \,;\, \tau^{-k}(x) \in [(k_{1,n}-i)\delta_{1,n},(k_{1,n}-i+1)\delta_{1,n}) \right\} \\
            &=&  \sharp \left\{ 1\leq k\leq n \,;\, x_0-\tau^{-k}(x) \in [t/n^3+(i-1)\delta_{1,n},t/n^3+i\delta_{1,n}) \right\},\\
    \mbox{and} \quad m_j^2(x)&:=& \sharp \left\{ 1\leq k\leq n \,;\, \tau^{-k}(x) \in [x_0+t/n^3+(j-1)\delta_{2,n},x_0+t/n^3+j\delta_{2,n}) \right\} .
   \end{array} \]
We then obtain
\[ \begin{array}{rcl}
  L_{n,x} &=& e^{n(X_0+X'_0)} \prod_{k ; \tau^{-k}(x) \in I_1}(x_0-\tau^{-k}(x))\prod_{k ; \tau^{-k}(x) \in I_2}(\tau^{-k}(x)-x_0) \\
            &\geq& e^{n(X_0+X'_0)} \prod_{i=1}^{k_{1,n}}\left(t/n^3+(i-1)\delta_{1,n}\right)^{m_i^1(x)} \prod_{j=1}^{k_{2,n}}\left(t/n^3+(j-1)\delta_{2,n}\right)^{m_j^2(x)}  \\
            &=& e^{n(X_0+X'_0)} \prod_{i=0}^{k_{1,n}-1}\left(t/n^3+i\delta_{1,n}\right)^{m_{i+1}^1(x)} \prod_{j=0}^{k_{2,n}-1}\left(t/n^3+j\delta_{2,n}\right)^{m_{j+1}^2(x)} .
        \end{array} \]
Using the discrepancy we get, for $i=1,\hdots,k_{1,n}$ and $j=1,\hdots,k_{2,n}$
\[ \frac{m_i^1(x)}{n}-\delta_{1,n} \leq D_n \quad \mbox{and} \quad \frac{m_j^2(x)}{n}-\delta_{2,n} \leq D_n .\]
Hence,
\[ L_{n,x} \geq e^{n(X_0+X'_0)} \left[\prod_{k=0}^{k_{1,n}-1}\left(\frac{t}{n^3}+k\delta_{1,n}\right)\right]^{n(D_n+\delta_{1,n})}  \left[\prod_{k=0}^{k_{2,n}-1}\left(\frac{t}{n^3}+k\delta_{2,n}\right)\right]^{n(D_n+\delta_{2,n})}=L_n .\]
Let us look at the first term of the product:
\[ \begin{array}{rcl}
     Q_1=\left[\prod_{k=0}^{k_{1,n}-1}\left(\frac{t}{n^3}+k\delta_{1,n}\right)\right]^{n(D_n+\delta_{1,n})} &=& \left[\delta_{1,n}^{k_{1,n}} \left(\frac{t}{n^3 \delta_{1,n}}\right)\prod_{k=1}^{k_{1,n}-1}\left(\frac{t}{n^3\delta_{1,n}}+k\right)\right]^{n(D_n+\delta_{1,n})}\\
      &=& \left[\frac{t \delta_{1,n}^{k_{1,n}-1}}{n^3} \frac{\Gamma\left(\frac{t}{n^3\delta_{1,n}}+k_{1,n}\right)}{\Gamma\left(\frac{t}{n^3\delta_{1,n}}+1\right)}\right]^{n(D_n+\delta_{1,n})}. 
\end{array} \]
Choose $\delta_{1,n}$ such that $\delta_{1,n}n^3 \rightarrow \infty$. Since for all $x\in[1,2]$ $\Gamma(x) \leq 1$, it comes
\[ Q_1 \geq \left(\frac{t\delta_{1,n}^{k_{1,n}-1}}{n^3}\right)^{n(D_n+\delta_{1,n})}\exp\left(n(D_n+\delta_{1,n})\ln\Gamma\left(\frac{t}{n^3\delta_{1,n}}+k_{1,n}\right)\right) .\]
But $\frac{t}{n^3\delta_{1,n}}+k_{1,n}=x_0/\delta_{1,n}$ and (\ref{dev_gamma}) give
\[ \begin{array}{rcl}
    \ln\Gamma\left(\frac{t}{n^3\delta_{1,n}}+k_{1,n}\right) &=& \ln\Gamma\left(\frac{x_0}{\delta_{1,n}}\right) \geq \ln\left(\left(\frac{x_0}{\delta_{1,n} e}\right)^{x_0/\delta_{1,n}} \sqrt{2\pi x_0/ \delta_{1,n}}\right) \\
      &=& -\frac{x_0}{\delta_{1,n}}\ln\delta_{1,n} -(1-\ln x_0)\frac{x_0}{\delta_{1,n}}+\frac{1}{2}\ln(2\pi x_0)-\frac{1}{2}\ln\delta_{1,n} .\\
    \end{array} \]
Hence we obtain
\[ \begin{array}{rcl}
 Q_1 &\geq&  \left(\frac{t\delta_{1,n}^{k_{1,n}-1}}{n^3}\right)^{n(D_n+\delta_{1,n})} \exp(n(D_n+\delta_{1,n})\left[-(x_0/\delta_n +1/2)\ln\delta_{1,n}\right. \\
     && \hspace*{9cm} \left. -X_0/\delta_{1,n}+\frac{1}{2}\ln(2\pi x_0)\right]) \\
     &=& \left(\frac{t}{n^3}\right)^{n(D_n+\delta_{1,n})} \exp([k_{1,n}-x_0/\delta_n-3/2]n(D_n+\delta_{1,n})\ln\delta_{1,n}-X_0n \\
     && \hspace*{4cm} - nX_0D_n/\delta_{1,n}+n/2(D_n+\delta_{1,n})\ln(2\pi x_0)) \\
     &=& \left(\frac{t}{n^3}\right)^{n(D_n+\delta_{1,n})} \exp\left(-X_0n-nX_0D_n/\delta_{1,n}+n/2(D_n+\delta_{1,n})\ln x_0\right) \\
     && \hspace*{3cm} \exp\left(-[\frac{t}{n^3\delta_{1,n}}+\frac{3}{2}]n(D_n+\delta_{1,n})\ln\delta_{1,n} +n/2(D_n+\delta_{1,n})\ln(2\pi)\right) \\
     &\geq& \left(\frac{t}{n^3}\right)^{n(D_n+\delta_{1,n})} \exp\left(-X_0n-nX_0D_n/\delta_{1,n}+n/2(D_n+\delta_{1,n})\ln x_0\right)  \\
     &=& \exp\left(n(D_n+\delta_{1,n})\ln\left(\frac{t}{n^3}\right)-nX_0-\frac{nX_0D_n}{\delta_{1,n}}+n/2(D_n+\delta_{1,n})\ln x_0\right) \\
     &\geq& \exp\left(n(D_n+\delta_{1,n})\ln t-nX_0-\frac{nX_0D_n}{\delta_{1,n}}-2n(D_n+\delta_{1,n})\ln n^3\right) \\
     &\geq& \exp\left(n(D_n+\Delta_n)\ln t-nX_0-\frac{nX_0D_n}{\delta_n}-2n(D_n+\Delta_n)\ln n^3\right).
 \end{array} \]
Doing the same with the second term of the product we get:
\[ \begin{array}{rcl}
 Q_2 &=& \left[\prod_{k=0}^{k_{2,n}-1}\left(\frac{t}{n^3}+k\delta_{2,n}\right)\right]^{n(D_n+\delta_{2,n})} \\
     &\geq& \exp\left(n(D_n+\Delta_n)\ln t-nX'_0-\frac{nX'_0D_n}{\delta_n}-2n(D_n+\Delta_n)\ln n^3\right) .
   \end{array}\]
Therefore it comes:
\[ L_n \geq \exp\left(2n(D_n+\Delta_n)\ln t-\frac{2nD_n}{\delta_n}-4n(D_n+\Delta_n)\ln n^3\right) ,\]
and finally we obtain, using also (\ref{upper_bound_Delta_n}):
\[ \begin{array}{rcl}
  \frac{1}{L_n} &\leq& \exp(-2n(D_n+\Delta_n)\ln t)\exp\left(4n(D_n+\Delta_n)\ln n^3 + \frac{2nD_n}{\delta_n}\right) \\
                &\leq& \exp(-2n(D_n+\Delta_n)\ln t)\exp\left(\frac{n}{2\ln^2 n}+4n(D_n+\delta_n)\ln n^3 + \frac{2nD_n}{\delta_n}\right) .
   \end{array} \]
Set $\widetilde{w_2}(n)=4n(D_n+\delta_n)\ln n^3 + \frac{2nD_n}{\delta_n}$, we want that $\widetilde{w_2}(n)/n^2$ is summable. This is the case if $\delta_n$ and $D_n$ can be chosen such that
\[ \begin{array}{c}
  \frac{D_n}{\delta_n} \leq \frac{1}{\ln^{1+\epsilon}n} \\
  (D_n+\delta_n)\ln n \leq \frac{1}{\ln^{1+\epsilon}n} ,
  \end{array} \]
where $\epsilon$ is an arbitrary positive constant.\\
Thus, if
\[ D_n = \frac{1}{\ln^{3+\epsilon}n} ,\]
and consequently we choose $\delta_n=1/\ln^{2+\epsilon_1} n$ with $\epsilon_1<\epsilon$, then $\widetilde{w_2}(n)/n^2$ is summable. Moreover, it comes:
\[ \frac{1}{L_n} \leq \exp(-3\frac{n}{\ln^{2+\epsilon_1} n}\ln t)\exp\left(\frac{n}{2\ln^2 n}+\frac{25n}{2\ln^{2+\epsilon_1}n}+\frac{2n}{\ln^{1+\epsilon-\epsilon_1}n}\right) .\]
Taking $\epsilon_1=\epsilon/2$ we get:
\[ \frac{1}{L_n} \leq \exp(-3\frac{n}{\ln^{2+\epsilon/2} n}\ln t)\exp\left(15\frac{n}{\ln^{1+\epsilon/2}n}\right) = C_{n,t}\exp(\widetilde{w}(n)) ,\]
with
\[ C_{n,t}=\exp\left(-\frac{1}{5}\widetilde{w}(n)\ln t\right) \]
(see (\ref{wtilde})). We can remark that the cases $x_0=0$ and $x_0=1$ are included in the previous work. 

Now assume that $v(x)=\exp( s_0(X_0+X_0') ) (\sigma_0 (x-x_0))^{s_0}$ with $\sigma_0 (x)=x$ or $\sigma_0 (x)= |x|$ and $s_0 >0$, i.e. we allow to have a positive power $s_0$ (note that the absolute value may be needed for this to make sense for certain powers $s_0$). Then working as above we obtain that if the discrepancy of $\tau$ satisfies condition (\ref{cond_disc}) then 
\[ \begin{array}{l}
   U_{n,x} \leq \exp(s_0\widetilde{w}(n)) \quad \mbox{for all}\ x\in [0,1) \quad\mbox{and} \\
   1/L_{n,x} \leq C_{n,t}^{s_0}\exp(s_0\widetilde{w}(n)) \quad \mbox{for all}\ x\in E_{x_0,t}.
   \end{array}\]
Now, assume more generally, that $v(x)= \exp\left( \sum_{i=0}^l s_i (X_i +X_i') \right) \prod_{i=0}^l ( \sigma_i (x-x_i))^{s_i}$ as in (\ref{P}) and (\ref{X}). We denote $S=\sum_{i=0}^l x_i$. Then working as above we get that for all $x\in E_t=\cap_{i=0}^{l}E_{x_i,t}$, 
\begin{equation}\label{lowerbound} 
\frac{1}{L_{n,x}} \leq C_{n,t}^{S}\exp(S\widetilde{w}(n)) ,
\end{equation}
and we still have that $|E_t| \xrightarrow[t\to0]{}1$.\\
Furthermore, for all $x\in[0,1)$ we also have 
\begin{equation}\label{upperbound} 
U_{n,x} \leq \exp(S\widetilde{w}(n)) .
\end{equation}
From now on let $w(n)=\frac{n}{\ln^{1+\epsilon/4}n}$. Then $w$ is subadditive, satisfies $\sum_{n\geq1}w(n)/n^2<\infty$, and $\widetilde{w}(n)/w(n)\to 0$. This last condition implies that, for $n$ large enough, say $n\geq n(t)$, we get for $x\in E_t$,
\begin{equation} \label{275}
 \frac{1}{L_{n,x}} \leq \exp(w(n)) ,
\end{equation}
and therefore, for $f\in L^p([0,1])$ it comes
\begin{equation}\label{upperboundnorm}
\left\|T^nf\right\|_{L^p(E_t)} \leq \exp(w(|n|))\|f\|_{L^p([0,1])} \quad \mbox{for}\ |n|\geq n(t).
\end{equation}
\endpf

\noindent\textbf{Definition of the functional calculus:} For $\phi\in A_w$ and $f\in L^p([0,1])$, using Theorem~\ref{bounds} we obtain
\[ \sum_{|n| \geq n(t)}\left|\widehat{\phi}(n)\right|\left\|T^nf\right\|_{L^p(E_t)} \leq \sum_{|n|\geq n(t)}\left|\widehat{\phi}(n)\right|e^{w(|n|)}\left\|f\right\|_{L^p([0,1])} < \infty, \]
which shows that the series $\sum_{n\in\Z}\widehat{\phi}(n)T^nf$ converges absolutely in $L^p(E_t)$. If $\phi^t(T)f$ denotes the infinite sum $\sum_{n \in \Z}\widehat{\phi}(n) T^nf$ in $L^p(E_t)$ then notice that for $0<t_2<t_1<1$ we have that $E_{t_1} \subset E_{t_2}$ and $(\phi^{t_2}(T)f)|_{E_{t_1}}= \phi^{t_1}(T)f$. Since $| \cup_{0<t<1}E_t| =1$, we denote by $\phi(T)f$ the function defined almost everywhere on $[0,1]$ such that $(\phi(T)f)|_{E_t} = \phi^t(T)f$ for all $0<t<1$.

We now give a second way of looking at $\phi(T)f$. Let $S_m=\sum_{|n| \leq m}\widehat{\phi}(n)T^nf$ be the partial sum of series $\sum_{n \in \Z} \widehat{\phi}(n)T^nf$ . We will prove that $S_m$ tends to $\phi(T)f$ in measure.
Let $(k_n)$ be an increasing subsequence of $\N$ and $t_n \in (0,1)$ be a decreasing sequence tending to $0$, and so $E_{t_n}\nearrow E$ with $|E|=1$. We proved previously that $S_m \rightarrow \phi(T)(f)$ in $L^p(E_t)$ for all $t\in(0,1)$, and so we can choose inductively increasing subsequences of positive integers $(k_n)_n \supseteq (k_n^1)_n \supseteq (k_n^2)_n \supseteq \hdots$ such that $(S_{k_n^i})_n$ converges almost everywhere on $E_{t_i}$ to the function $(\phi(T)f)|_{E_{t_i}}$ for all $i$. Thus the ``diagonal'' subsequence $(S_{k_n^n})_n$ converges almost everywhere to $\phi(T)f$. This proves that every subsequence $(S_{k_n})_n$ of $(S_n)$ has a further subsequence $(S_{k_n^n})_n$ wich converges $a.e.$ to the function $\phi(T)f$. Since the measure of $[0,1]$ is finite, this proves that $(S_n)$ converges in measure to $\phi(T)f$.\\[0.3cm]
\noindent\textbf{Properties of the functional calculus:} We will show that the above functional calculus has enough properties to prove the existence of an hyperinvariant subspace for the operator $T$. Namely we will prove the following:
\begin{prop}
The above functional calculus satisfies:
\begin{enumerate}[(I)]
\item There exists a non-empty set $D \subset L^p([0,1])\setminus \{0\}$ such that for all $\phi\in A_w$ and $f\in D$ we have $\phi(T)f \in L^p([0,1])$.
\item If $\phi_n\rightarrow \phi$ in $A_w$ and $f_n \rightarrow f$ in $L^p([0,1])$ then $\phi_n(T)f_n \rightarrow \phi(T)f$ in measure.
\item If $\phi_n\rightarrow \phi$ in $A_w$ and $f\in D$ then $\phi_n(T)f\rightarrow \phi(T)f$ in $L^p([0,1])$.
\item If $\psi\in A_w$, $f\in D$ and $S\in \{T\}'$ then $\psi(T)(Sf)=S(\psi(T)f)$.
\item If $\phi, \psi\in A_w$ and $f\in D$ then
\[ \phi(T)\left(\psi(T)f\right)=(\phi\psi)(T)f .\]
\item If $\psi\in A_w\setminus\{0\}$ and $f\in L^p([0,1])\setminus\{0\}$ satisfies $\psi(T)f\in L^p([0,1])$ then $\psi(T)f\neq 0$.
\end{enumerate}
\end{prop}

\beginpf  $(I)$ Let $G_t=E_t^c=\cup_{i=0}^l\cup_{n\geq1}E_{n,x_i,t}^c$, (see (\ref{setEnxt}) and (\ref{P})), we have $|G_t|\xrightarrow[t\to0]{}0$. We consider the two sets
\[ D_t=\{ f\in L^p([0,1]) ; f=0 \ \mbox{on} \ G_t\} \quad \mbox{and} \quad D=\cup_{t>0}D_t .\]
We remark that $D$ is dense in $L^p([0,1])$ for $p<\infty$.\\
Fix $f\in D$. There exists $t>0$ such that $f\in D_t$. We claim that
\begin{equation} \label{claim}
\left| T^{-m} f (x) \right| \leq \left( \sup_{y \in E_t} \frac{1}{L_{m,y}} \right)
\left| f(\tau^{-m} (x)) \right| \quad \mbox{ for all } m \in \N \mbox{ and }x \in [0,1].
\end{equation}
Indeed, if $x \in E_{m,x_i,t}^c$ for some $m \in \N$ and $i \in \{ 1, \ldots , l\}$ then there exists  $k\in\{1,\hdots,n\}$ such that $|\tau^{-k}(x)-x_i|<t/m^3$. But
\[ |\tau^{-k}(x)-x_i|=|\tau^{m-k}(\tau^{-m}(x))-x_i| < \frac{t}{m^3} < \frac{t}{(m-k)^3} ,\]
so $\tau^{-m}(x) \in E_{m-k,x_i,t}^c \subset G_t$. Thus $f(\tau^{-m}(x))=0$ for $x\in E_{m,x_i,t}^c$ and hence (\ref{claim}) is valid in this case. On the other hand, if $x \in \cap_{m \in \N} \cap_{i=1}^l E_{m,x_i,t} (=E_t)$ then (\ref{claim}) is valid by (\ref{upperboundnorm}) This finishes the proof of (\ref{claim}). 

Now (\ref{claim}) and (\ref{275}) imply that 
\begin{equation} \label{claim2}
\left| T^{-m}f (x)\right| \leq e^{w(m)}\left| f(\tau^{-m}(x) \right|
\quad \mbox{ for all } m > n(t) \mbox{ and }x \in [0,1].
\end{equation}
Combining (\ref{claim}), (\ref{claim2}) and Theorem~\ref{bounds} we obtain that 
for $f\in D_t$ and $\phi\in A_w$
\[ \begin{array}{rcl}
  \sum_{m\in\Z}\left|\widehat{\phi}(m)\right| \left\|T^mf\right\|_{L^p([0,1])} &\leq&  \sum_{m=0}^{n(t)} \left| \widehat{\phi}(m) \right| \left\| T^m  
                \right\| \| f \|_{L^p([0,1])} \\
                && + \sum_{m=-n(t)}^{-1}\left| \widehat{\phi}(m)\right| \sup_{x \in E_t}  
                \frac{1}{L_{n,x}}\|f\|_{L^p([0,1])} \\
                && + \sum_{|m| > n(t)}\left|\widehat{\phi}(m)\right| e^{w(m)}\left\|f\right\|_{L^p([0,1])} <\infty . 
\end{array}\]
Thus the series $\sum_{m\in\Z}\widehat{\phi}(m)T^mf$ converges absolutely in $L^p([0,1])$ for $f\in D_t$.\\[0.2cm]
\indent $(II)$ We have for $t\in(0,1)$
\[ \begin{array}{l}
  \left\|\sum_{|m|>n(t)}\widehat{\phi_n}(m)T^mf_n - \sum_{|m|>n(t)}\widehat{\phi}(m)T^mf\right\|_{L^p(E_t)}  \\  
        \hspace*{1cm} \leq  \sum_{|m|>n(t)}|\widehat{\phi_n}(m)-\widehat{\phi}(m)|\|T^mf_n\|_{L^p(E_t)} + \sum_{|m|>n(t)}|\widehat{\phi}(m)|\|T^m(f_n-f)\|_{L^p(E_t)} \\
        \hspace*{1cm} \leq \sum_{|m|>n(t)}|\widehat{\phi_n}(m)-\widehat{\psi}(m)|e^{w(|m|)}\sup_n\|f_n\|_{L^p(E_t)} \\
        \hspace*{6cm}  + \sum_{|m|>n(t)}|\widehat{\psi}(m)|e^{w(|m|)}\|f_n-f\|_{L^p([0,1])} \\
        \hspace*{1cm} \leq \|\phi_n-\phi\|_{A_w}\sup_n\|f_n\|_{L^p(E_t)}+\|\phi\|_{A_w}\|f_n-f\|_{L^p([0,1])} \xrightarrow[n\to\infty]{} 0
  \end{array} \]
Also,
\[ \sum_{m=0}^{n(t)}\left|\widehat{\phi}(m)\right|\|T^m(f_n-f)\|_{L^p(E_t)} \leq \sum_{m=0}^{m(t)}\left|\widehat{\phi}(m)\right|\|T^m\| \|f_n-f\|_{L^p([0,1])} \xrightarrow[n\to\infty]{}0.\]
Finally, for $x\in E_t$ we have
\[ \begin{array}{rcl}
 \sum_{m=1}^{n(t)}\left|\widehat{\phi}(-m)\right||T^{-m}(f_n-f)(x)| &\leq&  
         \sum_{m=1}^{n(t)}\left|\widehat{\phi}(-m)\right|\frac{1}{L_{m,x}}\left|(f_n-f)(\tau^{-m}x)\right| \\
   &\leq& \sum_{m=1}^{n(t)}\left|\widehat{\phi}(-m)\right|\sup_{x\in E_t}\frac{1}{L_{m,x}}\left|(f_n-f)(\tau^{-m}x)\right| ,
   \end{array}\]
which converges to 0 in measure as $n\to\infty$. Since $|E_t|\to 1$ as $t\to 0$, the above estimates imply $(II)$.\\[0.2cm]
\indent $(III)$ For $f\in D$ there exists $t >0$ such that $f\in D_t$. So
\[ \begin{array}{rcl}
  \|\phi_n(T)f-\phi(T)f\|_{L^p([0,1])} &\leq& \sum_{m\in\Z}\left|\widehat{(\phi_n-\phi)}(m)\right|\|T^mf\|_{L^p([0,1])} \\
           &\leq& 
\sum_{m=0}^{n(t)} \left| \widehat{(\phi_n-\phi)}(m) \right| \left\| T^m \right\| \| f \|_{L^p([0,1])} \\
& & + \sum_{m=-n(t)}^{-1}\left| \widehat{(\phi_n-\phi)}(m) \right| 
\sup_{x \in E_t} \frac{1}{L_{m,x}} \| f \|_{L^p([0,1])}\\                 
  & & +\sum_{|m| > n(t)}\left|\widehat{(\phi_n-\phi)}(m)\right|e^{w(|m|)}\|f\|_{L^p([0,1])}\\
  & & \mbox{ (by (\ref{claim}), (\ref{claim2}) and Theorem~\ref{bounds})} \\
           &\leq& \|\phi_n-\phi\|_{A_w}\|f\|_{L^p([0,1])} \xrightarrow[n\to\infty]{}0 ,
    \end{array}\]
which implies $(III)$.\\[0.2cm]
\indent $(IV)$ and $(V)$ These two equalities are true for trigonometric polynomials. Since they are dense in $A_w$, there exist $(\phi_n)$ and $(\psi_n)$ sequences of trigonometric polynomials such that $\phi_n \rightarrow\phi$ and $\psi_n\rightarrow\psi$ in $A_w$.\\
By $(II)$ with $f_n=f$ we have $\psi_n(T)(Sf) \rightarrow \psi(T)(Sf)$ in measure and furthermore, by $(III)$, we have also $S(\psi_n(T)f)\rightarrow S(\psi(T)f)$ in $L^p([0,1])$. So $(IV)$ is true.\\
Since $A_w$ is a Banach algebra, we have $\phi_n\psi_n\rightarrow \phi\psi$ in $A_w$. By $(II)$ and $(III)$ we get that $\phi_n(\psi_n(T)f)\rightarrow \phi(\psi(T)f)$ in measure, and also by $(III)$, we have $(\phi_n\psi_n)(T)f\rightarrow (\phi\psi)(T)f$ in $L^p([0,1])$. So $(V)$ is true.\\[0.2cm]
\indent $(VI)$ We take $\psi \geq 0$. Since $\psi\neq0$ there exists $t_0$ such that $\psi(e^{it_0})>0$, and by continuity of $\psi$, there exists an interval $I$ of length $\delta$ centered at $t_0$ such that $\psi(e^{it})>0$ for all $t$ on $I$.\\
Let $R_\delta$ be the rotation of angle $\delta$, and set $\psi_k=\psi\circ R_\delta^k$, \emph{i.e.} $\psi_k(e^{it})=\psi(e^{i(t+k\delta)})$. We recover the circle by successively shift $I$ by the rotation $R_\delta$. Then, given $N=\E(2\pi/\delta)+1$, we have $\sum_{k=0}^N \psi_k(e^{it})>0$ for all $t\in[0,2\pi)$. So we built a function on $A_w$, positive on $\T$, and in particular which does not vanish on $\T$. This function has an inverse $\Psi$ in $A_w$. We get:
\begin{equation}\label{eq}
\Psi(T)\sum_{k=0}^N \psi_k(T)f=f .
\end{equation}
Moreover, 
\[ \begin{array}{rcl}
\widehat{\psi_k}(n) &=& \frac{1}{2\pi}\int_0^{2\pi}\psi(e^{i(t+k\delta)})e^{-int}dt \\
      &=& \frac{1}{2\pi}\int_0^{2\pi}\psi(e^{i\theta})e^{-in\theta}e^{ik\delta}d\theta \\
      &=& e^{ik\delta}\widehat{\psi}(n),
\end{array}\]
hence
\[ \begin{array}{rcl}
\psi_k(T)f &=& \sum_{n\in\Z}\widehat{\psi_k}(n)T^nf \quad \mbox{(in measure)}\\
          &=& e^{ik\delta}\sum_{n\in\Z}\widehat{\psi}(n)T^nf \quad \mbox{(in measure)}\\
          &=& e^{ik\delta}\psi(T)f.
\end{array}\]
Thus the equality $(\ref{eq})$ becomes
\[ \left(\sum_{k=0}^Ne^{ik\delta}\right) \Psi(T)\psi(T)f=f ,\]
which implies, since $f$ is not the zero function, that $\psi(T)f \neq 0$.\\
\endpf
\vspace*{0.4cm}

\noindent \textbf{The hyperinvariant subspace:} Since the Beurling algebra is regular, there exist $\phi, \psi \in A_w$ with $\psi, \phi \geq 0$, such that $\phi\psi=0$. Set
\[ M=\left\{f\in L^p([0,1]) \;;\; \phi(T)(Sf)=0 \; \forall\, S\in \{T\}'\right\} ,\]
and $M_{hi}=\overline{M}^{L^p}$.\\
By construction $M_{hi}$ is closed and hyperinvariant. It remains to show that $M_{hi}$ is non trivial.\\
First $M_{hi}\neq L^p([0,1])$. Otherwise, take $f\in L^p([0,1])\setminus\{0\}$, there exists a sequence $f_n\in M$ such that $f_n\rightarrow f$ in $L^p$. Then by $(II)$ (and taking $S=Id$) we obtain $\phi(T)f_n\rightarrow \phi(T)f$ in measure and so $\phi(T)f=0$, which implies a contradiction by $(VI)$.\\
To prove that $M_{hi}\neq \{0\}$, we construct a non zero element in $M_{hi}$. By $(I)$ there exists $g\in D$ such that $\psi(T)g \in L^p([0,1])$. Then, by $(IV)$ and $(V)$, 
\[ \phi(T)(S\psi(T)g)=S(\phi(T)(\psi(T)g))=S((\phi\psi)(T)g)=0 .\]
This finishes the proof  Theorem~\ref{main}.\\
\endpf

\begin{rem}
For a function  $f$ defined on $[0,1]$ one may introduce the mixing-class associated to $f$ by
\[ \begin{array}{rcl}
 \mathcal{F}_f &:=& \left\{ \theta:[0,1]\rightarrow \R \;;\; \mbox{there exist a partition of}\ [0,1], \mbox{say} \cup_{i=0}^r[a_i,a_{i+1})\ \mbox{and}\right. \\
            && \left.\ c \ \mbox{a permutation of this partition such that}\ \theta (x)=f(c(x)) \right\} .
\end{array} \]
For example, for $\alpha \in[0,1]$, the function $\theta(x)=\{x+\alpha\}$ belongs to $\mathcal{F}_x$ where we consider the partition $[0,1]=[0,\alpha)\cup[\alpha,1)$. If $f$ is a bijection, $\theta$ is too. Our main Theorem~\ref{main} remains valid for more general weights $v$, namely if the functions $x-x_i$ ($1 \leq i \leq l$) in the definition of $v$ in (\ref{P}) are replaced by functions $\theta_i \in \mathcal{F}_{x-x_i}$. The spectral radius of the associated operator $T$ does not change in this case and all the previously estimates are obtained again by working on permuted subintervals of $[0,1]$. 
\end{rem}

\section{Example of ergodic transformation}
The most classical example of ergodic transform is the irrational rotation $\tau(x)=\{x+\alpha\}$, and we then obtain operator called Bishop-type operator. Davie \cite{Davie74} first proved that for almost all irrational number, namely the non Liouville numbers (they are dense in $\R$ with Lebesgue measure equals to 0), the Bishop operator associated to the weight $v(x)=x$ has hyperinvariant subspace. Later result due to MacDonald \cite{GMD90} generalized the result to a larger class of weight but for the same kind of irrational numbers as Davie. In \cite{Flattot}, the work of Davie was generalized to a larger class of irrationals, but for the weight $v(x)=x^s$, with $s$ a positive real. A recent work of Chalendar and Partington \cite{CP08} generalize the previous result for weight of the form $v(x)=\prod_{k=1}^K \{x-\beta_k\}^{\gamma_k}$ with $\gamma_k>0$ for some irrationals (including some Liouville numbers). This operator includes the product of Bishop-operator.

In this section we give sufficient conditions on an irrational number $\alpha$ such that the discrepancy of the ergodic transformation $\tau (x) = \{ x +\alpha \}$ satisfies condition (\ref{cond_disc}) and thus our Theorem~\ref{main} applies.  Recall that for a real number $t$ we denote by $<t>$ the distance from $t$ to the nearest integer, \emph{i.e.} $<t>= \min_{n \in {\mathbb Z}} |t-n|$. Also recall that if $\psi$ is a non-decreasing positive function defined on the positive integers, then an irrational number $\alpha$ is said to be of type $< \psi$ if $q<q \alpha > \geq 1/\psi(q)$ for all $q\in\N$. This is a measure of ``irrationality'' of the number $\alpha$. The smaller the function $\psi$ is, the ``farther away'' is $\alpha$ from the rationals; the larger the $\psi$ is, the ``closer'' $\alpha$ is allowed to be to the rationals.\\
Using results in \cite{KN74} we can prove the following:
\begin{prop} \label{psi}
Let $\alpha$ be an irrational number of type $< \psi$ where $\psi(q)= \exp(q^{1/(3+\epsilon)})$ for some $\epsilon >0$. Then
\[ D_N(\alpha) = O\left(\frac{1}{\ln^{3+\epsilon/3}N}\right) .\]
\end{prop}

\beginpf
In order to prove this proposition we first recall the two lemmas in \cite[p122-123]{KN74}:
\begin{enumerate}[(a)]
\item The discrepancy of $\omega=(n\alpha)$ satisfies
\[D_N(\omega) \leq C\left(\frac{1}{m}+\frac{1}{N}\sum_{h=1}^m\frac{1}{h<h\alpha>}\right) ,\]
for any positive integer $m$.
\item Let $\alpha$ be of type $<\psi$. Then,
\[ \sum_{h=1}^m\frac{1}{h<h\alpha>}=O\left(\psi(2m)\ln m+ \sum_{h=1}^m\frac{\psi(2h)\ln h}{h}\right) .\]
\end{enumerate}
We remark that the function $h \mapsto \frac{\psi(2h)\ln h}{h}$ is non-decreasing, so it comes
\[ \begin{array}{rcl}
   \sum_{h=1}^m\frac{\psi(2h)\ln h}{h} &\leq& \int_1^{m+1} \frac{\psi(2h)\ln h}{h}dh \\
          &=& O\left( \int_1^{m+1} e^{(2h)^{1/(3+\epsilon)}}\frac{2}{3+\epsilon}(2h)^{1/(3+\epsilon)-1} dh \right) \\
          &=& O\left( e^{(2m+2)^{1/(3+\epsilon)}} \right) .
    \end{array} \]
Thus 
\[ \begin{array}{rcl}
 \sum_{h=1}^m\frac{1}{h<h\alpha>} &=& O\left(e^{m^{1/(3+\epsilon/2)}}\right) + O\left( e^{(2m+2)^{1/(3+\epsilon)}} \right) \\
                      &=& O\left(e^{m^{1/(3+\epsilon/2)}}\right) ,
 \end{array} \]
and therefore we obtain
\[ D_N(\omega) = O\left(\frac{1}{m}+\frac{1}{N}e^{m^{1/(3+\epsilon/2)}}\right) \quad \mbox{for all positive integers} \ m.\]
Now choose $m= \lfloor \ln^{3+\varepsilon/2}(N\ln^{-3-\epsilon/3}N) \rfloor$ to obtain that $D_N(\alpha) = O\left( \frac{1}{\ln^{3+\varepsilon/3}N}\right)$ and finish the proof of the proposition.\\
\endpf

We recall the definition:
\begin{defn}
The irrational number $\alpha$ is a Liouville number if and only if $\alpha$ is not of type $< \phi$ for any power function $\phi$, \emph{i.e.} of the
form $\phi (q)=q^n$ (where n is a fixed positive integer). Equivalently, if for all integer $n$, there exist some integers $p$ and $q$ with $q > 1$ satisfying
$0<\left|\alpha-p/q\right|<1/q^n$.
\end{defn}

Looking closer at the definition of type $<\psi$, we remark that if $\alpha$ is of type $<\psi$ then $|a-p/q|\geq 1/(\psi(q)q^2)$ for all $p\in\N$. Thus we can use the results in \cite{Flattot} to estimate the size of the set $A$ of the irrational numbers $\alpha$ of type $< \psi$ where $\psi (q)= \exp ( q^{1/(3+\epsilon)})$ for some $\epsilon >0$. We have
\begin{prop}
For $f(x)=1/\ln^8(x/2)$ the $f$-Hausdorff measure of $A^c$ is zero.\\
Furthermore, as soon as $g$ converges faster than $f$ to $0$ in $0$, we have $H^g(A^c)=0$.
\end{prop}

And we can explicit Liouville number for which the operator has hyperinvariant subspace.
\begin{prop}
Let $b\geq2$ be an integer, and let $(u_n)$ be a sequence of positive integers satisfying, for $n$ large enough, the two conditions
\[ \begin{array}{l}
  nu_n +\frac{\ln \beta}{\ln b}< u_{n+1} \quad \mbox{with} \quad \beta=\frac{b}{b-1} ,\\
  u_{n+1} < \frac{b^{u_n/(3+\epsilon)}}{\ln b}.
   \end{array}\]
Then the number $\alpha=\sum_{n\geq0}1/b^{u_n}$ is a Liouville number of type $< \psi$ where $\psi (q)= \exp ( q^{1/(3+\epsilon)})$. For this choice of $\alpha$ we obtain by Proposition~\ref{psi} and Theorem~\ref{main} that $T$ has a nontrivial hyperinvariant subspace.
\end{prop}

\noindent \textbf{Example:} Taking $b=10$ and $u_n=n!$, we obtain the classical example of Liouville number, and for this one the operator $T$ has a nontrivial hyperinvariant subspace.

\section{Case higher dimension}
The previously described approach  allows us to easily extend Theorem~\ref{main} to the case of weighted composition operators
on $L^p([0,1]^d)$ for $d$ a positive integer. We consider the operator 
\[ \begin{array}{rcl}
 T :  L^p([0,1]^d)& \rightarrow & L^p([0,1]^d)\\
  f & \mapsto & v\, f\circ \tau
\end{array} \]
where $v\in L^{\infty}([0,1]^d)$ and $\tau$ is a bijective ergodic transformation on $[0,1]^d$. We assume that $v$ and $\tau$ can be written as
\[ \begin{array}{l}
  v(x)=\prod_{i=1}^d v_i(x_i)  \quad \mbox{and}\\
 \tau(x)=(\tau_1(x_1),\hdots,\tau_d(x_d)) ,
 \end{array}\]
where $v_i \in \mathcal{P}$ for all $i\in\{1,\hdots,d\}$.
Such an example $\tau$ of ergodic transformation is the rotation on the $d-$dimensional torus with angle $\alpha=(\alpha_1,\ldots,\alpha_d)$, defined by
\[ \begin{array}{rcl} 
\tau :   [0,1]^d& \rightarrow & [0,1]^d\\
  x=(x_1,\ldots,x_d)& \mapsto & (\{x_1+\alpha_1\},\ldots,\{x_d+\alpha_d\}). 
\end{array} \]
A vector $\alpha=(\alpha_1,\ldots,\alpha_d)$ is said irrational if $1,\alpha_1,\ldots,\alpha_d$ are linearly independent over $\mathbb{Q}$, and if $\alpha$ is irrational, then the rotation $\tau$ is known to be uniquely ergodic so in particular ergodic. For this ergodic transform it was proved in \cite{CFG08} that the operator has no eigenvalues for all irrational $\alpha$, but it was not known about the existence of hyperinvariant subspaces. 
Also, MacDonald extends in \cite{GMD91} its first work to obtain 
hyperinvariant subspace for some Bishop-type operator with a non-vanishing weight $v$ and an irrational rotation $\tau$.

In order to apply the previous work, we remark that the only change is in the bounds of $T^n$, but considering these ergodic transformation and weight allows to easily obtain that
\[ \begin{array}{l}
  \left\|T^nf\right\|_{L^p([0,1]^d)} \leq e^{dw(n)}\left\|f\right\|_{L^p([0,1]^d)} \quad \mbox{for } n\geq 0,\\
  \left\|T^nf\right\|_{L^p(E_t)} \leq C_{n,t}^{d} e^{dw(n)}\left\|f\right\|_{L^p([0,1]^d)} \quad \mbox{for } n< 0 ,
  \end{array}\]
where $E_t$ is defined as before but $E_{n,x_0,t}$ becomes
\[ E_{n,x_0,t}=\{ x\in[0,1]^d \,;\, |\tau^{-k}_i(x)-x_0| \geq t/n^3 \; \forall k=1,\hdots,n,i=1\hdots d \} ,\]
and assuming that 
\[ \sup_{1\leq i\leq d} D_{i,n} =O\left(\frac{1}{\ln^{3+\epsilon}n}\right) 
\quad \mbox{for some } \epsilon>0,\]
where $D_{i,n}$ denotes the discrepancy associated to $\tau_i$.
Under these assumptions we obtain to the existence of a hyperinvariant subspace for this operator.

\begin{rem}
\begin{enumerate}
\item Theorem~\ref{main} should remain valid for more general weights $v$ where the functions $x-x_i$ in (\ref{P}) are replaced by any function whose graph is the union of linear segments. However, some of the above calculations may become more technical when somebody tries to complete this task.
\item In the higher dimensional case, our main Theorem~\ref{main} (and the scheme of our proof) should remain valid when one considers more general maps $\tau:[0,1]^d \to [0,1]^d$ where the discrepancy is defined by replacing $[\alpha , \beta)$ in (\ref{DN}) by $\prod_{k=1}^d [\alpha_i, \beta_i)$. In that case one should need to create grids of $[0,1]^d$ using $d$-dimensional cubes in order to obtain the bounds of $T^n$ for $n \in \Z$. However, this may make the proof and the notation more tedious. 
\end{enumerate}
\end{rem}

\end{document}